\documentclass[final,3p,times]{elsarticle}
\usepackage{amsmath,amssymb,amsfonts,epsf,amsthm}
\usepackage{hyperref}
\hypersetup{colorlinks=true, urlcolor=blue, citecolor=cyan, pdfborder={0 0 0},}
\usepackage{soul}
\usepackage{natbib}
\usepackage{multirow}

\usepackage[ruled, vlined, Algorithm]{algorithm}
\usepackage{algorithmic}

\newtheorem{theorem}{Theorem}[section]
\newtheorem{lemma}{Lemma}[section]
\newtheorem{proposition}{Proposition}[section]

\newtheorem{corollary}{Corollary}[section]

\newtheorem{definition}{Definition}[section]

\numberwithin{equation}{section}

\usepackage{amssymb}

\journal{--}

\allowdisplaybreaks[1]

\begin{document}

\begin{frontmatter}

\title{  Some algebraic  properties of a  class of integral graphs  determined by their spectrum }

\author[label1]{Jia-Bao Liu}
\ead{liujiabaoad@163.com;liujiabao@ahjzu.edu.cn}
\author[label2]{S.Morteza Mirafzal}
\ead{Morteza$\_$mirafzal@yahoo.com}
\author[label3]{Ali Zafari \corref{1}}
\ead{zafari.math.pu@gmail.com}
\address[label1]{School of Mathematics and Physics,
Anhui Jianzhu University, Hefei 230601, P.R. China}
\address[label2]{Lorestan University, Department of Mathematics,  Faculty of Science,
Khorramabad, Iran}
\address[label3]{Department of Mathematics, Faculty of Science,
Payame Noor University, P.O. Box 19395-4697, Tehran, Iran}
\cortext[1]{Corresponding author}
\begin{abstract}
Let $\Gamma=(V,E)$ be a graph. If all the eigenvalues of the adjacency matrix of the graph $\Gamma$ are integers, then we say that $\Gamma$ is an integral graph. A graph $\Gamma$ is determined by its spectrum if every graph cospectral to it is in fact isomorphic to it. In this paper, we investigate some algebraic properties of the Cayley graph $\Gamma=Cay(\mathbb{Z}_{n}, S)$,  where $n=p^m$, ($p$ is a  prime  integer, $m\in\mathbb{N}$) and $S=\{{a}\in\mathbb{Z}_{n}\,|\,\, (a, n)=1\}$. First, we show that $\Gamma$ is an integral graph. Also we determine the automorphism group of $\Gamma$. Moreover, we show that  $\Gamma$ and $K_v \bigtriangledown\Gamma$ are
determined by their spectrum.
\end{abstract}
\begin{keyword}
Adjacency spectrum, Laplacian spectrum, Integral graph, DS graph.
\MSC[2010] 05C50, 05C31
\end{keyword}
\end{frontmatter}
\section{Introduction}
\label{sec:introduction}
The graphs in this paper are simple, undirected and connected. We always assume that
$\overline{\Gamma}$ denotes the complement graph of $\Gamma$.
The eigenvalues of a graph $\Gamma$ are the eigenvalues of the adjacency matrix of $\Gamma$.
The spectrum of $\Gamma$ is the list of the eigenvalues of the adjacency matrix
of $\Gamma$ together with their multiplicities, and it is denoted by Spec$(\Gamma)$, see ~\cite{book3}.
If all the eigenvalues of the adjacency matrix of the graph $\Gamma$ are integers,
then we say that $\Gamma$ is an integral graph. The notion of integral graphs was
first introduced by F. Harary and A. J. Schwenk in 1974, see ~\cite{paperf.1}.  In general, the problem of
 characterizing integral graphs seems to be very difficult. There are good surveys in this area, see ~\cite{paperb-1}. For more results depend on the integral graphs  and their applications in engineering networks see ~\cite{pap.s.b,pap.a,pap.s.m}.
For any vertex $v$ of a connected graph $\Gamma$,
  we denote the set of vertices of $\Gamma$ at distance $r$ from $\Gamma$  by $\Gamma_r(v)$. Then we have,
$$\Gamma_r(v)=\{u\in V(\Gamma)\,\  | \,\,d(u, v)=r\},$$ where $d(u, v)$ denotes the distance in $\Gamma$ between the vertices $u$ and $v$, and
 $r$ is a non negative integer not exceeding $d$, the diameter of $\Gamma$.
It is clear that $\Gamma_0(v)=\{v\}$, and $V(\Gamma)$ is partitioned  into the disjoint subsets $\Gamma_0(v), ..., \Gamma_d(v)$, for each $v$ in $V(\Gamma)$. The graph $\Gamma$ is called distance regular with diameter $d$ and intersection array $\{b_0, ..., b_{d-1}; c_1, ..., c_d\}$ if it is regular of valency $k$ and, for any two vertices $u$ and $v$ in $\Gamma$ at distance $r$  we have $|\Gamma_{r+1}(v)\cap\Gamma_{1}(u)|=b_r$, $(0\leq r\leq d-1)$  and $|\Gamma_{r-1}(v)\cap\Gamma_{1}(u)|=c_r$  $(1\leq r\leq d)$. The intersection numbers $c_r, b_r$ and $a_r$ satisfy
$a_r=k-b_r-c_r\,\ \,\,  \,\, (0\leq r\leq d),$ where $a_r$ is the number of neighbours of $u$ in $\Gamma_r(v)$.
Let $G$ be a finite group and $H$ a subset of $G$  such that it is closed under taking inverses and does not contain the identity.
A Cayley graph $\Gamma=Cay(G, H)$ is the graph whose vertex set and edge set are defined as follows:
$$V (\Gamma) =G; \,\, \,\ E(\Gamma) = \{\{x, y\} \, | \,\, x^{-1}y \in H\}.$$
It is well known that
if $\Gamma$ is a distance regular graph with  valency  $k$, diameter $d$, adjacency matrix $A$, and intersection array
$$\{b_0, b_1, ..., b_{d-1}; c_1, c_2, ..., c_d\}.$$
then, the tridiagonal $(d+1)\times(d+1)$ matrix,
\begin{center}
$	B=\begin{bmatrix}
a_0 & b_0 &0 &0 &...& \\
c_1 & a_1 & b_1&0&...& \\
0& c_2 & a_2&b_2& & \\
& & &...& \\
& & & & \\
& & & &c_{d-2}&a_{d-2}&b_{d-2}&0 \\
& & &...&0&c_{d-1}&a_{d-1}&b_{d-1} \\
& & &...&0&0&c_d&a_d \\
\end{bmatrix},$
\end{center}
determines all the eigenvalues of $\Gamma$  ~\cite{book2-1}.
Note that, the concept of distance regular graphs date back to the 1960s. They were
defined by Biggs, see  ~\cite{book1-1}, and their basic theory was developed by him and others.
Distance regular graphs of diameter $2$ are just the connected strongly regular
graphs. The theory of distance regular
graphs has connections to many parts of graph theory such as, design theory, coding theory,
geometry, and group theory.
Two graphs with the same spectrum are called cospectral. It is not hard to see that the spectrum of a graph
does not determine its isomorphism class. The authors in ~\cite{paper2} proposed the question: which
graphs are determined by their spectrum? It seems hard to prove a graph to be determined by its spectrum. Up to now, only
a few classes of  graphs are proved to be determined by their spectrum, such as: the path $P_n$, the complete graph $K_n$ and the cycle $C_n$,
graph $Z_n$ and their complements, see   ~\cite{paperzz11,paperzz14,paperzz16}.
For a graph $\Gamma$, let $A(\Gamma)$ and $L(\Gamma)=D(\Gamma) - A(\Gamma)$
be respectively the adjacency matrix and Laplacian matrix of $\Gamma$, where $D(\Gamma)$ is the diagonal matrix of vertex degrees with $\{d_1, d_2, ..., d_n\}$ as diagonal entries.
 Laplacian spectrum and their applications are involved indiverse
theoretical problems on complex networks ~\cite{k-1,z-1}. Many results
have been devoted to studying Laplacian spectrum for complex
networks ~\cite{j-1,pap.j.b.2}. Calculating the Laplacian spectrum  of networks
has many applications in lots of aspects, such as the topological
structures and dynamical processes ~\cite{c-1}.
Algebraic properties of various classes of Cayley graphs have been studied by various authors, see  ~\cite{paper6,paper7}.
 In this paper, we want to study some algebraic properties of a class of Cayley graphs constructed on the cyclic additive group $\mathbb{Z}_{n}$, denoted by $\Gamma=Cay(\mathbb{Z}_{n}, S)$,  where $n=p^m$, ($p$ is a  prime  integer, $m\in\mathbb{N}$) and $S=\{{a}\in\mathbb{Z}_{n}\,|\,\, (a, n)=1\}$. It is easy to check that $S$ is an inverse closed subset in the group $\mathbb{Z}_{n}$ and $0 \notin S$. Thus, $\Gamma$ is a simple graph.  This class of graphs is a espacial subclass of graphs, which are investigated from some other aspects by
 Basi\'{c} and  Ili\'{c}  ~\cite{paper3a}.   Using the theory of distance
regular graphs,  we show that the adjacency  spectrum  of  $\Gamma$ is $\{ n-p^{m-1}, 0^{(n-p)}, (-p^{m-1})^{(p-1)}\}$, where the superscripts give the multiplicities of eigenvalues with multiplicity greater than one. Finally, we show that any graph cospectral with the multicone graph $K_v\bigtriangledown \Gamma$ is determined by their adjacency
spectrum as well as their Laplacian spectrum, where $K_v$ is the complete graph on $v$ vertices.
\section{Definitions And Preliminaries}
\begin{definition}\label{b.1}
~\cite{book2-1,book1}
Let $\Gamma$ be a graph with automorphism group $Aut(\Gamma)$.
We say that $\Gamma$ is a vertex transitive graph if for all vertices $ x, y$ of $\Gamma$, there is an automorphism $\theta$ in $Aut(\Gamma)$ satisfying $\theta(x)=y$.  Also,
 we say that $\Gamma$ is distance transitive graph if, for all vertices $u, v, x, y$ of $\Gamma$ such that $d(u, v)=d(x, y)$, there is an automorphism $\theta$ in $Aut(\Gamma)$ satisfying $\theta(u)=x$ and $\theta(v)=y$
\end{definition}
\begin{theorem} \label{b.2}
~\cite{l.b}
 Let $\Gamma$ be a graph such that contains $k$ components $\Gamma_{1}, ... , \Gamma_{k}$. If  for any $i\in I=\{1, ... , k\}$, we have  $\Gamma_{i}\cong \Gamma_{1}$ then $Aut(\Gamma)\cong Aut(\Gamma_{1})wr_{I} Sym(k)$, where the wreath product is defined.
\end{theorem}
\begin{definition} \label{b.3}
~\cite{paper8}
Let $\Gamma_1\cup\Gamma_2$ denote the disjoint union of graphs $\Gamma_1$ and $\Gamma_2$. The join  $\Gamma_1\bigtriangledown\Gamma_2$ is the graph obtained from $\Gamma_1\cup\Gamma_2$ by joining every vertex of $\Gamma_1$ with every vertex of $\Gamma_2$. A multicone graph is defined to be the join of a clique and a regular graph.
\end{definition}
\begin{theorem} \label{b.4}
~\cite{paper2}
If $\Gamma$ is a distance regular graph with diameter $d$ and girth $g$ satisfying
one of the following properties, then every graph cospectral with $\Gamma$ is also
distance regular, with the same parameters as $\Gamma$:
\begin{itemize}
\item[(i)]
 $g \geq 2d - 1$,
\item[(ii)]
 $g \geq 2d - 2$ and $\Gamma$ is bipartite.
\end{itemize}
\end{theorem}
\begin{proposition}\label{b.5}
~\cite{paper2}
For regular graphs, being $DS$ (or not $DS$) is equivalent for the adjacency matrix, the adjacency matrix of the complement, and the Laplacian matrix.
\end{proposition}
\begin{proposition}\label{b.6}
~\cite{paper2}
The following graph and its complement, which have at most
four eigenvalues, are regular $DS$ graphs:
\begin{itemize}
\item[(i)]
The disjoint union of $k$ copies of a strongly regular $DS$ graph.
\end{itemize}
\end{proposition}
\begin{theorem} \label{b.9}
~\cite{paper4}
Let $\Gamma_1$ and $\Gamma_2$ be two graphs with the Laplacian spectrum $\lambda_1\geq\lambda_2\geq...\geq\lambda_n$ and $\mu_1\geq\mu_2\geq...\geq\mu_m$, respectively. Then, the Laplacian spectrum
of  $\Gamma_1\bigtriangledown\Gamma_2$,  is
$n+m, m+\lambda_1, m+\lambda_2,..., m+\lambda_{n-1}, n+\mu_1, n+\mu_2, ..., n+\mu_{m-1}, 0$.
 \end{theorem}
\begin{theorem} \label{b.10}
~\cite{book3}
Let $\Gamma$ be a graph on $n$ vertices. Then, $n$ is a Laplacian
eigenvalue of $\Gamma$ if and only if $\Gamma$ is the join of two graphs.
 \end{theorem}
\begin{lemma}\label{b.11}
~\cite{book3}
A connected graph $\Gamma$ has exactly one positive eigenvalue if and only if it is a complete
multipartite graph.
\end{lemma}
\section{Main results}
\begin{theorem}\label{c.1}
Let $\Gamma=Cay(\mathbb{Z}_{n} , S)$ be the Cayley graph on the cyclic  group $\mathbb{Z}_{n}$, where $n=p^m$, ($p$ is a  prime  integer and $m\in\mathbb{N}$)  and $S=\{{a}\in\mathbb{Z}_{n}\,|\,\, \text {(a, n)=1}\}$.  Then,
 $$Aut(\Gamma)\cong Sym(p^{m-1})wr_{I} Sym(p),$$  where $I=\{1, 2, ... , p\}.$
\end{theorem}
\begin{proof}
Let $V(\Gamma)=\{1, ... , n\}$ be the vertex set of $\Gamma$. Note that if $m = 1$, then the result immediately follows. Because in this case, $\Gamma\cong K_p$, where $K_p$  is the complete graph
on $p$ vertices. Hence in the sequel, we assume that $m\geq2$. Let $T=\langle p \rangle = \{kp \  | \  0 \leq k \leq p^{m-1}-1 \}$  be the subgroup of the group $\mathbb{Z}_{n}$ of order $p^{m-1}$. It is clear that $T$ and every co-set of $T$ is an independent set in the graph $\Gamma$.
In fact if  $T+a$ is a co-set of $T$ in the group $\mathbb{Z}_{n}$  such that
   $T \cap T+a= \emptyset $, then $a$ and $p$ are co-prime   and hence we have $a\in S$.
This follows that every co-set of $T$ is a clique of order $p^{m-1}$ in the complement of the graph $\Gamma$.
 Thus, $\overline{\Gamma}$  contains of $p$ disjoint components $\Gamma_{1},\Gamma_{2}, ... , \Gamma_{p}$ such that  $\Gamma_{i}\cong K_{p^{m-1}}$ ($1\leq i\leq p$), where $K_{p^{m-1}}$ is the complete graph on $p^{m-1}$ vertices.  This follows that  $\overline{\Gamma}\cong pK_{p^{m-1}}$. Hence, by Theorem \ref{b.2}, $Aut(\overline{\Gamma})\cong Aut(K_{p^{m-1}}) wr_{I} Sym{(p)}= Sym(p^{m-1}) wr_{I} Sym{(p)}$.
 On the other hand, it is well known that for any graph $\Gamma$, $Aut(\Gamma)=Aut(\overline{\Gamma})$, see ~\cite{book3}.
\end{proof}
\begin{proposition}\label{c.2}
Let $\Gamma=Cay(\mathbb{Z}_{n}, S)$ be the Cayley graph on the cyclic  group $\mathbb{Z}_{n}$, where $n=p^m$, ($p$ is a  prime  integer and $m\in\mathbb{N}$)  and $S=\{{a}\in\mathbb{Z}_{n}\,|\,\, \text {(a, n)=1}\}$.  Then $\Gamma$ is a  distance transitive graph.
\begin{proof}
 Suppose $u, v, x, y$ are vertices of $\Gamma$ such that $d(u, v)=d(x, y)=r$, where $r$ is a non negative integer not exceeding $d$, the diameter of  $\Gamma$. So $d(u, v)=d(x, y)=1$ or 2, since we now have the diameter of  $\Gamma$ is $d=2$. In the following cases we show that $\Gamma$ is a  distance transitive graph.\newline

Case 1.
If $d(u, v)=d(x, y)=2$, then $u^{-1}v\notin S$ and $x^{-1}y\notin S$.  Therefore,  two vertices $u$ and $v$ are adjacent in  the complement $\overline{\Gamma}$ of $\Gamma$, also two vertices $x$ and $y$ are adjacent in  the complement $\overline{\Gamma}$ of $\Gamma$. By  Theorem \ref{c.1}, we know that $\overline{\Gamma}$  contains $p$ components
$\Gamma_{1}, \Gamma_{2}, ... , \Gamma_{p}$ such that for any $i\in \{1, 2, ... , p\}$, $\Gamma_{i}\cong  K_{p^{m-1}}$. Therefore,
 $\overline{\Gamma}\cong pK_{p^{m-1}}$. If $u=x$, then $u, v, y$ lie in a clique of graph $\overline{\Gamma}$, and hence we may assume
$\theta=(vy) \in Aut(\overline{\Gamma})=Aut({\Gamma})$, so $\theta(u)=x$ and $\theta(v)=y$.
If  $u\neq x$ and $v\neq y$, then $u, v$ lie in a clique of graph $\overline{\Gamma}$, say $\Gamma_i$, also
$x, y$ lie in a clique of graph $\overline{\Gamma}$, say $\Gamma_j$, where  $\Gamma_i\neq\Gamma_j$ or $\Gamma_i=\Gamma_j$. Hence we may assume $\theta=(ux)(vy) \in Aut(\overline{\Gamma})=Aut({\Gamma})$. Thus $\theta(u)=x$ and $\theta(v)=y$.\newline

Case 2.
If $d(u, v)=d(x, y)=1$, then we can show that there is an automorphism $\theta$ in  $Aut(\Gamma)$ such that $\theta(u)=x$ and $\theta(v)=y$.
\end{proof}
\end{proposition}
\begin{proposition}\label{c.3}
Let $\Gamma=Cay(\mathbb{Z}_{n}, S)$ be the Cayley graph on the cyclic  group $\mathbb{Z}_{n}$, where $n=p^m$, ($p$ is a  prime  integer and $m\in\mathbb{N}$)  and $S=\{{a}\in\mathbb{Z}_{n}\,|\,\, \text {(a, n)=1}\}$. Then $\Gamma$ is an  integral graph.
\begin{proof}
It is well known that if $\Gamma$ is a distance transitive graph, then  $\Gamma$ is also  distance regular, see ~\cite{book1}. Now,
let $V(\Gamma)=\{1, 2, ... ,n\}$ be the  vertex set of $\Gamma$.  Consider the vertex  $v=n$ in $V(\Gamma)$, then $\Gamma_0(v)=\{n\}$,
 $\Gamma_1(v)=\{a\in V(\Gamma) \,|\, \, (a, n)=1\}$, and
$\Gamma_2(v)=\{a\in V(\Gamma) \,|\, \, (a, n)\neq 1\}$.
 Let   $u$ be the vertex in $V(\Gamma)$ such that  $d(u, v)=0$,  then $u=v=n$ and $|\Gamma_{1}(v)\cap\Gamma_{1}(u)|=n-p^{m-1}$. Hence $b_0=n-p^{m-1}$, and  by Definition  of distance regularity of graph we have $a_0=(n-p^{m-1})-b_0=0$. Also, if $u$ in $V(\Gamma)$ and  $d(u, v)=1$,
then two vertices $u, v$ are adjacent in $\Gamma$, so  $|\Gamma_{0}(v)\cap\Gamma_{1}(u)|=1$, and $|\Gamma_{2}(v)\cap\Gamma_{1}(u)|=p^{m-1}-1$. Hence $c_1=1$,   $b_1=p^{m-1}-1$ and $a_1=(n-p^{m-1})-b_1-c_1=n-2p^{m-1}$. Finally, if  $u$ in $V(\Gamma)$ and  $d(u, v)=2$, then two vertices $u, v$ are not adjacent in $\Gamma$, so
 $|\Gamma_{1}(v)\cap\Gamma_{1}(u)|=n-p^{m-1}$, hence  $c_2=n-p^{m-1}$ and  $a_2=(n-p^{m-1})-(n-p^{m-1})=0.$ Thus the intersection array of  $\Gamma$  is $\{ n-p^{m-1}, p^{m-1}-1; 1, n-p^{m-1}\}.$ Therefore,  the tridiagonal $(3)\times(3)$ matrix,
 \begin{center}
 $\begin{bmatrix}
a_0 & b_0 &0  \\
c_1 & a_1 & b_1 \\
 0& c_2 & a_2 \\
\end{bmatrix}=\begin{bmatrix}
0& n-p^{m-1} &0  \\
1 & n-2p^{m-1} & p^{m-1}-1 \\
 0& n-p^{m-1} & 0 \\
\end{bmatrix}$
  \end{center}
  determines all the eigenvalues of $\Gamma$. It is clear that all the eigenvalues of $\Gamma$ are
   $n-p^{m-1}, 0, -p^{m-1}$, and  their multiplicities  are $1, n-p, p-1$, respectively. Thus $\Gamma$ is an  integral graph.
 \end{proof}
 \end{proposition}
 \begin{corollary}\label{c.4}
 Let $\Gamma=Cay(\mathbb{Z}_{n}, S)$ be the Cayley graph on the cyclic  group $\mathbb{Z}_{n}$, where $n=p^m$, ($p$ is a  prime  integer and $m\in\mathbb{N}$)  and $S=\{{a}\in\mathbb{Z}_{n}\,|\,\, \text {(a, n)=1}\}$. Then  the adjacency  spectrum  of  $\Gamma$ is $\{ n-p^{m-1}, 0^{(n-p)}, (-p^{m-1})^{(p-1)}\}$.
 \end{corollary}
\begin{theorem}\label{c.5}
 Let $\Gamma=Cay(\mathbb{Z}_{n}, S)$ be the Cayley graph on the cyclic  group $\mathbb{Z}_{n}$, where $n=p^m$, ($p$ is a  prime  integer and $m\in\mathbb{N}$)  and $S=\{{a}\in\mathbb{Z}_{n}\,|\,\, \text {(a, n)=1}\}$. Then $\Gamma$ is a $DS$ graph
with respect to its  the  adjacency  spectrum.
\begin{proof}
We know that if $p$ is even prime  integer, then
$\Gamma$ is isomorphic to the bipartite graph $K_{p^{m-1}, p^{m-1}}$, and hence
the result immediately follows.

Now, let  $p$ is an odd prime integer, then  $\Gamma$ is not bipartite graph. In particular,
 $g \geq 2d - 1$, because  the diameter of $\Gamma$ is $2$,  and  the girth of  $\Gamma$ is $3$. Hence by Theorem \ref{b.4},
  every graph cospectral with $\Gamma$ is also distance regular, with the same parameters as $\Gamma$. Because by Proposition
    \ref{c.2}, we know that $\Gamma$ is a distance regular graph.  Thus,
 $\Gamma$ is a $DS$ graph with respect to its  the  adjacency  spectra. Because by Proposition \ref{b.6}, $\overline{\Gamma}$  contains disjoint union of $p$ copies of  the strongly regular $DS$ graph $K_{p^{m-1}}$, in addition
 the graph  $\Gamma$ and its complement, which have at most four eigenvalues.
 \end{proof}
 \end{theorem}
 \begin{proposition}\label{c.6}
Let $\Pi$ be a graph cospectral with the multicone graph $K_v \bigtriangledown\Gamma$  with respect to its adjacency matrix spectrum,  where $\Gamma=Cay(\mathbb{Z}_{n}, S)$  which is defined as before.
 Then $\Pi$ is a bidegreed graph. Also,
$$Spec(\Pi)=\{  0^{(n-p)}, {(-p^{m-1})}^{(p-1)}, -1^{(v-1)}, (\frac{M+\sqrt{M^2+4N}}{2}), (\frac{M-\sqrt{M^2+4N}}{2})\},$$
where $M=v-1+p^m-p^{m-1}$ and $N=p^m+p^{m-1}v-p^{m-1}$.
\begin{proof}
We can deduce the following from Theorem 2.1.8 in ~\cite{paper1} and Theorem 2.1 in ~\cite{paper3}.
 \end{proof}
\end{proposition}
\begin{theorem}\label{c.8}
Consider the multicone graph $K_v \bigtriangledown\Gamma$, where $\Gamma=Cay(\mathbb{Z}_{n}, S)$
 which is defined as before.  Then $K_v \bigtriangledown\Gamma$ is $DS$ with respect to its adjacency matrix spectrum.
\begin{proof}
In the following, we proceed by induction on the number of vertices in $K_v$.  Let $K_v$ have one vertex and
$\Pi$ be  a graph cospectral with  the multicone graph $K_1 \bigtriangledown\Gamma$ with respect to its adjacency matrix spectrum.
By Proposition \ref{c.6}, it is easy to see that $\Pi$ has one vertex of degree $p^m$, say $j$.
Hence, if $Spec(\Pi-j)=Spec(\Gamma)$, then $\Pi-j \cong \Gamma$.
 Because, by Theorem \ref{c.5}, we know that  $\Gamma$  is  $DS$ graph
with respect to its  the  adjacency matrix spectrum.
So $\Pi\cong K_1 \bigtriangledown\Gamma$.
We assume inductively that this claim holds for $K_v$
that is, if $\Pi_1$ is a graph cospectral with the multicone graph $K_v \bigtriangledown\Gamma$ with
respect to its adjacency matrix spectrum, then  $\Pi_1\cong K_v \bigtriangledown\Gamma$.
 We show that the claim is true for $K_{v+1}$
that is,  if $\Pi$ is a graph  cospectral with the multicone graph
 $K_{v+1} \bigtriangledown\Gamma$ with respect to its adjacency matrix spectrum,  then $\Pi\cong K_{v+1} \bigtriangledown\Gamma$.
It is obvious that $\Pi$ has one vertex and $p^m+v$ edges more than $\Pi_1$. On the other hand,
by Proposition \ref{c.6}, we know that $\Pi_1$ has $v$ vertices of degree $p^m+v-1$ and $p^m$ vertices of
 degree $p^m-p^{m-1}+v$, and also $\Pi$ has $v+1$ vertices of degree $p^m+v$ and $p^m$ vertices
 of degree $p^m-p^{m-1}+v+1$. So, we must have $\Pi\cong K_1 \bigtriangledown\Pi_1$.
Now, by assume induction  we  conclude that $\Pi\cong K_{v+1} \bigtriangledown\Gamma$, and complete the proof.
 \end{proof}
\end{theorem}
\begin{theorem}\label{c.9}
Consider the  complement $\overline{K_v \bigtriangledown\Gamma}$  of multicone graph $K_v \bigtriangledown\Gamma$ with respect to its  adjacency spectrum, where $\Gamma=Cay(\mathbb{Z}_{n}, S)$  which is defined as before. Then $\overline{K_v \bigtriangledown\Gamma}$ is a $DS$ graph.
\begin{proof}
By Theorem \ref{c.1}, we know that
 $\overline{\Gamma}$  contains $p$ components $\Gamma_{1},\Gamma_{2}, ... , \Gamma_{p}$ such that
  $\Gamma_{i}\cong K_{p^{m-1}}$ ($1\leq i\leq p$). So  $\overline{\Gamma}\cong pK_{p^{m-1}}$. In addition, the adjacency matrix spectrum  of  $\overline{\Gamma}$ is $$\lbrace  (p^{m-1}-1)^{(p)}, -1^ {(p^{m}-p)}\rbrace.$$ Also, the adjacency matrix spectrum  of  $\overline{K_v}$ is $\{0^{(v)}\}$.
Thus the adjacency matrix spectrum  of  $\overline{\Gamma}\cup \overline{K_v}$ is $$\{  (p^{m-1}-1)^{(p)}, -1^ {(p^{m}-p)}, 0^{(v)}\}.$$
On the other hand, it  is not hard to see that   $\overline{\Gamma}\cup \overline{K_v}\cong\overline{K_v \bigtriangledown\Gamma}$.
Let $\Pi$ be a graph cospectral with the  complement $\overline{K_v \bigtriangledown\Gamma} $  of multicone graph $K_v \bigtriangledown\Gamma$ with respect to its  adjacency spectrum,
then $$Spec(\Pi)=Spec(\overline{K_v \bigtriangledown\Gamma})=\{  (p^{m-1}-1)^{(p)}, -1^ {(p^{m}-p)}, 0^{(v)}\}.$$
 It is easy to prove that $\Pi$ cannot be
regular, since regularity of a graph can be determined by its spectrum.
Also, we show that $\Pi$ is disconnected graph. Suppose to the contrary   $\Pi$ is connected,  hence by Lemma \ref{b.11},
$\Pi$ is complete multipartite graph, contradicting the adjacency spectrum of $\Pi$. Thus $\Pi$ is disconnected graph. Therefore, we conclude that $\overline{K_v \bigtriangledown\Gamma}$ is $DS$ with respect to its
  adjacency spectrum.
\end{proof}
\end{theorem}
\begin{proposition}\label{c.10}
Consider the multicone graph $K_v \bigtriangledown\Gamma$, where $\Gamma=Cay(\mathbb{Z}_{n}, S)$  which is defined as before.  Then $K_v \bigtriangledown\Gamma$ is $DS$ with respect to its  Laplacian spectrum.
\begin{proof}
 By  Theorem \ref{b.9}, the Laplacian matrix spectrum  of $K_v \bigtriangledown\Gamma$ is
$$\{ (n+v)^{(p+v-1)}, (n+v-p^{m-1})^{(n-p)},  0\}.$$
We proceed by induction on the number of vertices in $K_v$.  If $v= 1$, there is nothing to
prove.
We assume inductively that this claim holds for $K_v$
that is, if $Spec(L(\Pi_1))=
Spec(L(K_v \bigtriangledown\Gamma))$, then  $\Pi_1\cong K_v \bigtriangledown\Gamma$,
 where $\Pi_1$ is a graph cospectral with the multicone graph $K_v \bigtriangledown\Gamma$ with respect to its  Laplacian spectrum. We show that the claim is true for $K_{v+1}$
that is, if
$$Spec(L(\Pi))=
 Spec(L(K_{v+1} \bigtriangledown\Gamma))$$ $$=\{ (n+v+1)^{(p+v)}, (n+v+1-p^{m-1})^{(n-p)}, 0\},$$
 then $\Pi\cong K_{v+1} \bigtriangledown\Gamma$,
 where $\Pi$ is a graph cospectral with the multicone graph $K_{v+1} \bigtriangledown\Gamma$ with respect to its  Laplacian spectrum.
By Theorem \ref{b.10}, we know that  $\Pi_1$ and $\Pi$ are join of two graphs, because $n+v$ and $n+v+1$ are
eigenvalues of $\Pi_1$ and $\Pi$,  respectively. In addition,  $\Pi$ has one vertex of degree $n+v$ more than $\Pi_1$, say $j$, hence $Spec(L(\Pi-j))\cong Spec(L(K_v \bigtriangledown\Gamma))$, and by assume induction
$\Pi-j\cong K_v \bigtriangledown\Gamma$.
Thus it can be concluded $\Pi\cong  K_{v+1} \bigtriangledown\Gamma$.
 \end{proof}
\end{proposition}
\section{Conclusion}
In this paper, we computed the adjacency spectrum of a class of integral graphs, denoted by $\Gamma=Cay(\mathbb{Z}_{n}, S)$,  where $n=p^m$, ($p$ is a  prime  integer, $m\in\mathbb{N}$) and $S=\{{a}\in\mathbb{Z}_{n}\,|\,\, (a, n)=1\}$. Indeed, by  using the theory of distance
regular graphs,  it is shown that the adjacency  spectrum  of  $\Gamma$ is $\{ n-p^{m-1}, 0^{(n-p)}, (-p^{m-1})^{(p-1)}\}$, where the superscripts give the multiplicities of eigenvalues with multiplicity greater than one. Moreover, it is shown that the Cayley graph $\Gamma$ and $K_v \bigtriangledown\Gamma$ are
determined by their spectrum. Note that, this class of graphs is a espacial subclass of integral circulants, and hence  clearly this class not
only they are mathematically applicable, but also they are used in the design of engineering networks.

\bigskip
{\footnotesize
\noindent \textbf{Data Availability}\\
No data were used to support this study.\\[2mm]
\noindent \textbf{Conflicts of Interest}\\
The authors declare that there are no conflicts of interest
regarding the publication of this paper.\\[2mm]
\noindent \textbf{Acknowledgements}\\
This work was supported in part by Anhui Provincial Natural Science Foundation under Grant 2008085J01 and Natural Science Fund of Education Department of Anhui Province under Grant KJ2020A0478. \\[2mm]
{\footnotesize

\bigskip

\begin{thebibliography}{5}
\bibitem{book3}{C. Godsil and G. Royle}, \emph{Algebraic graph theory}, vol. 207 of Graduate Texts in Mathematics, Springer, New York, NY, USA, 2001.
\bibitem{paperf.1}{F.  Harary and  A. Schwenk}, \emph{Which graphs have integral spectra?}, in Graphs and Combinatorics, vol. 406 of Lecture Notes in Mathematics, pp. 45–51, Springer, 1974.
\bibitem{paperb-1}{K. Balinska, D. Cvetkovi\'{c}, Z. Rodosavljevi\'{c}, S. Simi\'{c}, and D. Stevanovi\'{c}},  A survey on integral graphs, Univ. Beograd, Publ. Elektrotehn. Fak. Ser. Mat, vol. 13, pp. 42–65, 2002.
 \bibitem{pap.s.b}{S. Blackburn and I. Shparlinski}, \emph{ On the average energy of circulant graphs}, Linear Algebra Appl, vol.428, pp.1956-1953, 2008.

\bibitem{pap.a}{ A. Ili\'{c} and M.  Basi\'{c}}, \emph{New results on the energy of integral circulant graphs}, Appl Math Comput, vol.218, pp.3470-3482, 2011.
  \bibitem{pap.s.m}{S. M. Mirafzal}, \emph{A new class of integral graphs constructed from the hypercube}, Linear Algebra and its Applications, vol.558,  pp.186-194, 2018.


 \bibitem{book2-1}{ A. E. Brower, A. M. Cohen, and A. Neumaier}, Distance-regular graphs, vol. 18, Springer, Berlin, Germany, 1989.
\bibitem{book1-1}{N. L. Biggs}, \emph{Intersection matrices for linear graphs},  Combinatorial Mathematics
and its applications, Proc. Oxford, 7-10 July,  pp.15-23, 1969.
\bibitem{paper2}{E. R. van Dam and  W. H. Haemers}, \emph{Which graphs are determined by their spectrum?}, Linear Algebra and its Applications,
vol.373, pp.241-272, 2003.
\bibitem{paperzz11}{X. G. Liu,  Y. P. Zhang, and X. Q. Gui}, \emph{The multi-fan graphs are determined by their Laplacian spectra}, Discrete Mathematics, vol.308,  pp.4267-4271, 2008.
\bibitem{paperzz14}{X. L. Shen, Y. P. Hou, and Y. P. Zhang}, \emph{Graph $Z_n$ and some graphs related to $Z_n$ are determined by their spectrum},
Linear Algebra and its Applications, vol.404,   pp.58-68, 2005.
\bibitem{paperzz16}{W. Wang and C. X. Xu}, \emph{The T-shape tree is determined by its Laplacian spectrum}, Linear Algebra and its Applications,
 vol.419,  pp.78-81, 2006.
\bibitem{k-1}{B. Y. Hou, H. J. Zhang, and L. Liu}, \emph{Applications of Laplacian spectra  for extended Koch
networks}, The European Physical Journal B, vol.85,30373, 2012.
\bibitem{z-1} {Z. Zhang, H. Liu, B. Wu, and T. Zou}, \emph{Spanning trees in  a fractal scale-free lattice},
Phys.Rev.E, vol.83, 016116, pp.1-8, 2011.
\bibitem{j-1}{J. B. Liu and X. F. Pan}, \emph{Asymptotic incidence energy of lattices}, Physica A: Statistical Mechanics and its Applications, vol.422,
pp.193-202, 2015.

\bibitem{pap.j.b.2}{J. B. Liu, Z. Y. Shi, Y. H. Pan, J. Cao, M. A. Aty, and  U. A. Juboori}, \emph{Computing the Laplacian spectrum of linear octagonal-quadrilateral networks and its applications}, Polycyclic Aromatic Compounds, Taylor and Francis Group, pp.1-12, 2020.
\bibitem{c-1}{Z. Cheng, J. Cao, and T. Hayat}, \emph{Cascade of failures in interdependent networks with
different average degree}, International Journal of Modern Physics C, vol.25, no 05, 1440006, 2014.
\bibitem{paper6}{S. M. Mirafzal and A. Zafari}, \emph{On the spectrum of a class of distance-transitive
graphs}, Electronic Journal of Graph Theory and Applications vol.5,  pp.63-69, 2017.
\bibitem{paper7}{ S. M. Mirafzal and A. Zafari}, \emph{An interesting property of a class of circulant graphs}, Journal of Mathematics, vol.2017, pp.1-4, 2017.
\bibitem{paper3a}{M.  Basi\'{c} and A. Ili\'{c}}, \emph{On the automorphism group of integral circulant graphs},
\emph{the electronic journal of combinatorics}, vol.18, pp.1-21,  2011.
\bibitem{book1}{ N. L. Biggs}, \emph{Algebraic graph theory}, New York, NY: Cambridge University Press, Cambridge, 1993.
\bibitem{l.b} {L. Beineke and  R. J Wilson}, \emph{Topics in algebraic graph theory}, Mathematical Sciences Faculty Publications,
\emph{Cambridge University Press}, 2004.
\bibitem{paper8} {J. F. Wang, H. Zhao, and Q. Haung}, \emph{Spectral charactrization of multicone graphs}, Czec. Math. J, vol.62,
pp.117-126, 2012.
\bibitem{paper4}{R. Merris}, \emph{Laplacian graph eigenvectors}, Linear Algebra Appl, vol.278,  pp.221-236, 1998.
\bibitem{paper1}{D. Cvetkovi\'{c}, P. Rowlinson, and S. Simi\'{c}},  \emph{An introduction to the theory of
graph spectra}, London Mathematical Society Student Texts, \textbf{75}. Cambridge
University Press, cambridge, 2010.
\bibitem{paper3}{Y. Hong, J. Shu, and K. Fang}, \emph{A sharp upper bound of the spectral radius of graphs},
J. Combinatorial Theory Ser. B, vol.81,  pp.177-183, 2001.


\end{thebibliography}
\end{document}